\documentclass[12pt]{amsart}
\usepackage{graphicx}
\usepackage{times}
\usepackage{mathptm}
\usepackage{eucal}
\newtheorem{thm}{Theorem}[section]  
\newtheorem{cor}[thm]{Corollary}
\newtheorem{quest}[thm]{Question}
\newtheorem{defin}[thm]{Definition} 
\newtheorem{lemma}[thm]{Lemma} 

\newtheorem{prop}[thm]{Proposition} 
 
\newcommand{\aaa}{\mbox{$\alpha$}}

\newcommand{\bbb}{\mbox{$\beta$}}
  
\newcommand{\Sss}{\mbox{$\Sigma$}}  
 
\newcommand{\rrr}{\mbox{$\rho$}} 
\newcommand{\Ggg}{\mbox{$\Gamma$}}
\newcommand{\Lll}{\mbox{$\Lambda$}}
\newcommand{\ggg}{\mbox{$\gamma$}}
 
\newcommand{\bdd}{\mbox{$\partial$}}

\newcommand{\Ddd}{\mbox{$\Delta$}}

\newcommand{\lll}{\mbox{$\lambda$}}

\begin{document}  



\title[Distance and bicompressibility]{Proximity in the curve 
complex: boundary reduction and bicompressible surfaces}   

\author{Martin Scharlemann}
\address{\hskip-\parindent
        Martin Scharlemann\\
        Mathematics Department\\
        University of California\\
        Santa Barbara, CA USA}
\email{mgscharl@math.ucsb.edu}

\thanks{Research partially supported by an NSF grant.} 

\date{\today}

\begin{abstract} 
Suppose $N$ is a compressible boundary component of a compact
irreducible
orientable $3$-manifold $M$ and $(Q, \bdd Q) \subset (M, \bdd M)$ is
an orientable properly embedded essential surface in $M$ in which some
essential component is incident to $N$ and no component is a disk.
Let $\mathcal{V}$ and $\mathcal{Q}$ denote respectively the sets of
vertices in the curve complex for $N$ represented by boundaries of
compressing disks and by boundary components of $Q$.  \bigskip

{\bf Theorem:} Suppose $Q$ is essential in $M$, then
$d(\mathcal{V}, \mathcal{Q}) \leq 1 - \chi(Q)$.
\bigskip

Hartshorn showed (\cite{Ha}) that an incompressible surface in a
closed $3$-manifold puts a limit on the distance of any Heegaard
splitting.  An augmented version of the theorem above leads to a
version of Hartshorn's result for merely compact $3$-manifolds.

\bigskip
In a similar spirit, here is the main result:

\bigskip

{\bf Theorem:} Suppose a properly embedded connected surface $Q$ is
incident to $N$.  Suppose further that $Q$ is separating and
compresses on both its sides, but not by way of disjoint disks.  Then 
either
    \begin{itemize} 
	\item $d(\mathcal{V}, \mathcal{Q}) \leq 1 - \chi(Q)$ or \item
	$Q$ is obtained from two nested connected incompressible
	boundary-parallel surfaces by a vertical tubing.
\end{itemize}

\bigskip

Forthcoming work with M. Tomova (\cite{STo}) will show how an
augmented version of this theorem leads to the same conclusion as in
Hartshorn's theorem, not from an essential surface but from an
alternate Heegaard surface.  That is, if $Q$ is a Heegaard splitting
of a compact $M$ then no other Heegaard splitting has distance greater
than twice the genus of $Q$.
\end{abstract}

\maketitle

\section{Introduction}

Suppose $N$ is a compressible boundary component of an orientable
irreducible $3$-manifold $M$ and $(Q, \bdd Q) \subset (M, \bdd M)$ is
an essential orientable surface in $M$ in which an essential component
is incident to $N$
and no component of $Q$ is a disk.  Let $\mathcal{V}, \mathcal{Q}$ denote
sets of vertices in the curve complex for $N$ represented respectively
by boundaries of compressing disks and by boundary components of $Q$.
We will show:

\begin{itemize}
\item The distance $d(\mathcal{V}, \mathcal{Q})$ in
the curve complex of $N$ is no greater than $1 - \chi(Q)$.
Furthermore, if no component of $Q$ is an annulus $\bdd$-parallel into
$N$, then
for each component $q$ of $Q \cap N$, $d(q, \mathcal{V}) \leq 1 -
\chi(Q)$.
\end{itemize}

A direct consequence is this generalization of a theorem of Hartshorn
\cite{Ha}:

\begin{itemize}  
   \item If $S$ is a Heegaard splitting surface for a compact
   orientable manifold $M$ and $(Q, \bdd Q) \subset (M, \bdd M)$ is a
   properly embedded incompressible surface, then $d(S) \leq 2 -
   \chi(Q)$.  
\end{itemize}

Both results are unsurprising, and perhaps well-known (see eg \cite{BS} 
for discussion of this in the broader setting of knots in bridge 
position with respect to a Heegaard surface).  

It would be of interest to be able to prove the second result 
(Hartshorn's theorem) for $Q$ a Heegaard surface, rather than an 
incompressible surface.  Of course this is hopeless in general: a second copy 
of $P$ could be used for $Q$ and that would in general provide no information 
about the distance of the splitting $P$ at all.  But suppose it is 
stipulated that $Q$ is not isotopic to $P$.  One possibility is that 
$Q$ is weakly reducible.  In that case (cf \cite{CG}) it is either the 
stabilization of a lower genus Heegaard splitting (which we revert 
to) or it gives rise to a lower genus incompressible surface and this 
allows the direct application of Hartshorn's theorem.  So in trying 
to extend Hartshorn's theorem to $Q$ a Heegaard surface, it suffices 
to consider the case in which $Q$ is strongly irreducible.  

The first step in extending \cite{Ha} to $Q$ a Heegaard surface is
carried out here, analogous in the program to the first result above.
Specifically, we establish that bicompressible but weakly
incompressible surfaces typically do not have boundaries that are
distant in the curve complex from curves that compress in $M$.

\begin{itemize}
   \item Suppose a properly embedded surface $Q$ is connected,
   separating and incident to $N$.  Suppose further that $Q$
   compresses on both its sides, but not by way of disjoint disks,
   then either
   \begin{itemize} 
       \item	$d(\mathcal{V}, \mathcal{Q}) \leq 1 - \chi(Q)$ or
       \item $Q$ is obtained from two nested connected boundary-parallel 
       surfaces by a vertical tubing. 
\end{itemize}
\end{itemize}

From this result forthcoming work will demonstrate, via a
two-parameter argument much as in \cite{RS}, that the genus of an
alternate Heegaard splitting $Q$ does indeed establish a bound on the
distance of $P$.

Maggy Tomova has provided valuable input to this proof.  Beyond
sharpening the foundational proposition (Propositions
\ref{prop:essential} and Theorem \ref{thm:bicompressible}) in a very
useful way, she provided an improved proof of Theorem
\ref{thm:Hartshorn}.

\section{Preliminaries and first steps}

First we recall some definitions and elementary results, most of 
which are well-known.
    
     \begin{defin} \label{defin:boundcomp}
	A $\bdd$-compressing disk for $Q$ is a disk $D \subset M$ so
	that $\bdd D$ is the end-point union of two arcs, $\aaa = D
	\cap \bdd M$ and $\bbb = D \cap Q$, and $\bbb$ is essential in
	$Q$. 
    \end{defin}

    \begin{defin} \label{defin:essential}
	A surface $(Q, \bdd Q) \subset (M, \bdd M)$ is
	{\em essential} if it is incompressible and some component is not
	boundary parallel.  An essential surface is {\em strictly essential}
	if it has at most one non-annulus component.  
	
	\end{defin}

	\begin{lemma} \label{lemma:boundredessential}
	     Suppose $(Q, \bdd Q) \subset (M, \bdd M)$ is a properly
	     embedded surface and $Q'$ is the result of
	     $\bdd$-compressing $Q$.  Then
	     \begin{enumerate}
	\item If $Q$ is incompressible so is $Q'$.
	\item If $Q$ is essential, so is $Q'$.

	\end{enumerate}	 
		    
   \end{lemma}
   
\begin{figure}[tbh]
   \centering
   \includegraphics[width=0.8\textwidth]{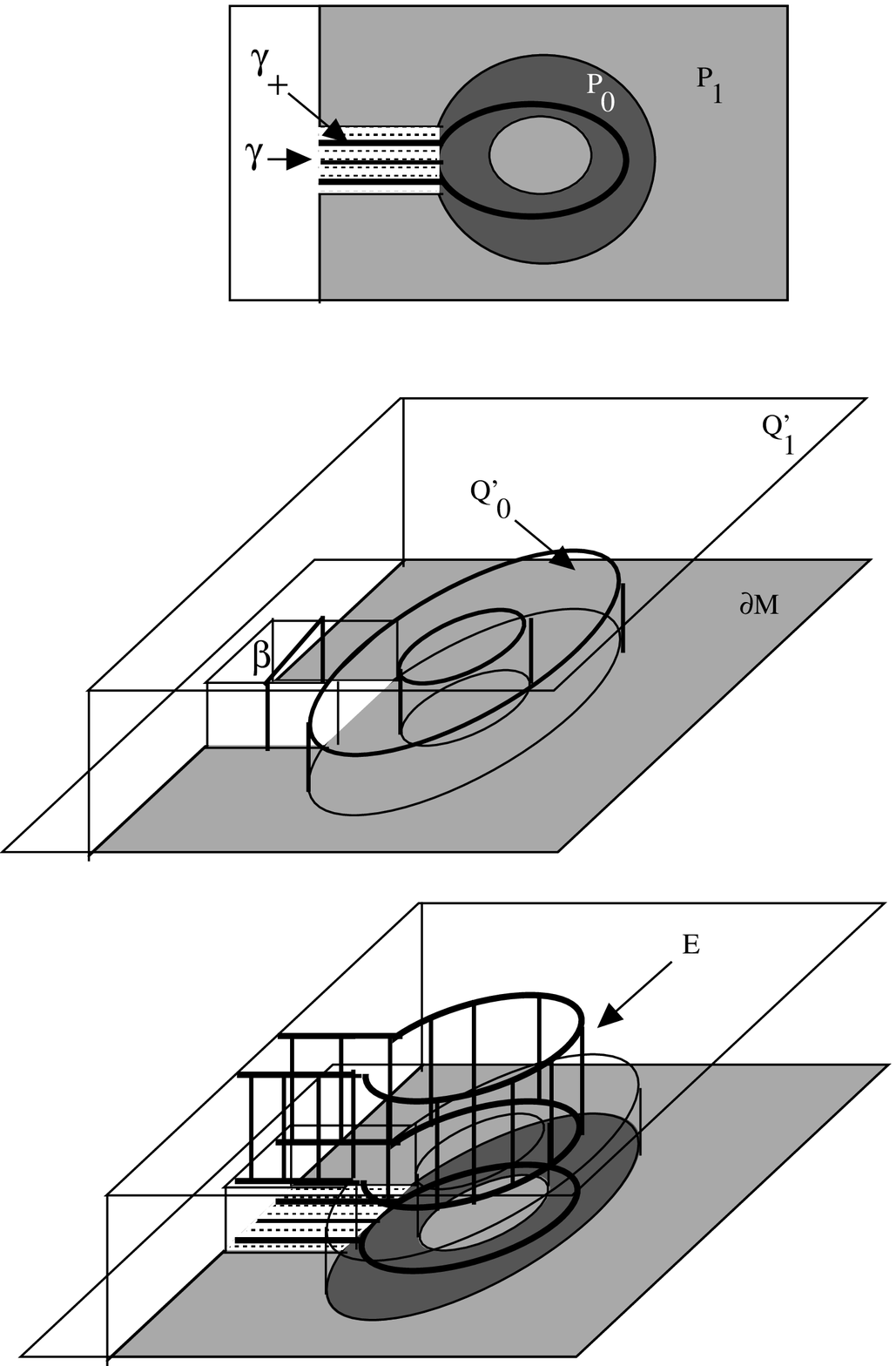}
   \caption{} \label{fig:essential}
   \end{figure}

   \begin{proof} A description dual to the boundary compression from
   $Q$ to $Q'$ is this: $Q$ is obtained from $Q'$ by tunneling along
   an arc $\ggg$ dual to the $\bdd$-compression disk. (The precise 
   definition of tunneling is given in Section \ref{sect:exam}.)
   Certainly any compressing disk for $Q'$ in $M$ is unaffected by
   this operation near the boundary.  Since $Q$ is incompressible, so
   then is $Q'$.  This proves the first claim.
     
   Suppose now that every component of $Q'$ is boundary parallel and
   the arc $\ggg$ dual to the $\bdd$-compression has ends on
   components $Q'_{0}, Q'_{1}$ of $Q'$ (possibly $Q'_{0} = Q'_{1}$). 
   If $\ggg$ is disjoint from the subsurfaces $P_{0}$ and $P_{1}$ of
   $\bdd M$ to which $Q'_{0}$ and $Q'_{1}$ respectively are parallel
   then tunneling along $\ggg$ merely creates a component that is
   again boundary parallel (to the band-sum of the $P_{i}$ along
   $\ggg$), contradicting the assumption that not all components of
   $Q$ are boundary parallel.  So suppose $\ggg$ lies in $P_{0}$, say. 
   If both ends of $\ggg$ lie on $Q'_{0}$ (so $Q'_{1} = Q'_{0}$) then
   the disk $\ggg \times I$ in the product region between $Q'_{0}$ and
   $P_{0}$ would be a compressing disk for $Q$, which contradicts the
   incompressibility of $Q$.  
   
   Finally, suppose $Q'_{1} \neq Q'_{0}$, so $P_{0} \subset P_{1}$ and
   $\ggg$ is an arc in $P_{1} - P_{0}$ connecting $\bdd P_{0}$ to
   $\bdd P_{1}$.  $P_{0}$ is not a disk, else the arc $\bbb$ in which
   the $\bdd$-compressing disk intersects $Q$ would not have been
   essential in $Q$.  So there is an essential simple closed curve
   $\ggg_{0} \subset P_{0}$ based at the point $\ggg \cap P_{0}$.
   Attach a band to $\ggg_{0}$ along $\ggg$ to get an arc $\ggg_{+}
   \subset P_{1}$ with both ends on $\bdd P_{1}$.  Then the disk
   $E_{1} = \ggg_{+} \times I$ lying between $P_{1} \subset \bdd M$
   and $Q'_{1}$ intersects $Q$ in a single arc, parallel in $M$ to $\ggg_{+}$
   and lying in the union of the top of the tunnel and $Q'_{0}$.  This
   arc divides $E_{1}$ into two disks; let $E$ be the one not incident
   to $\bdd M$.  $E$ then has its boundary entirely in $Q$ and since
   it is essential there, $E$ is a compressing disk for $Q$, again a
   contradiction.  See Figure \ref{fig:essential}.  From these various
   contradictions we conclude that at least one of the components of
   $Q'$ to which the ends of $\ggg$ is attached is not
   $\bdd$-parallel, so $Q'$ is essential.
   
   \end{proof}
   
\begin{defin} \label{defin:distance}
    
    Suppose $S$ is a closed orientable surface and $\aaa_{0},\ldots,
    \aaa_{n}$ is a sequence of essential simple closed curves in $S$ so
    that for each $1 \leq i \leq n$, $\aaa_{i-1}$ and $\aaa_{i}$ can
    be isotoped to be disjoint.  Then we say that the sequence is a
    length $n$ path in the curve complex of $S$ (cf \cite{He}).  
    
    The distance $d(\aaa, \bbb)$ between a pair $\aaa, \bbb$ of
    essential simple closed curves in $S$ is the
    smallest $n \in \mathbb{N}$ so that there is a path in the
    curve-complex from $\aaa$ to $\bbb$ of length $n$.  Curves are
    isotopic if and only if they have distance $0$.  
    
    Two sets of
    curves $\mathcal{V}, \mathcal{W}$ in $S$ have distance
    $d(\mathcal{V}, \mathcal{W}) = n$ if $n$ is the smallest distance
    from a curve in $\mathcal{V}$ to a curve in $\mathcal{W}$.
    
    \end{defin}
    
    \begin{prop} 
    \label{prop:essential}

	 Suppose $M$ is an irreducible compact
	 orientable $3$-manifold, $N$ is a compressible component of
	 $\bdd M$ and $(Q, \bdd Q) \subset (M, \bdd M)$ is a properly 
	 embedded essential surface with $\chi(Q) \leq 1$ and at least
	 one essential component incident to $N$.
    Let $\mathcal{V}$ be the set of essential curves in $N$ that bound 
    disks in $M$ and let $q$ be any component of $\bdd Q$.  
	\begin{itemize}
    \item If $Q$ contains an essential disk incident to $N$, then $d(\mathcal{V}, q) \leq 1$.
    \item If $Q$ does not contain any disk components, then $d(\mathcal{V}, q) \leq
	 1-\chi(Q)$ or $Q$ is strictly essential and $q$ lies in the boundary
	 of a $\bdd$-parallel
	 annulus component of $Q$.

    \end{itemize}
    \end{prop}

    \begin{proof}
     If $Q$ contains an essential disk $D$ incident to $N$, then $\bdd
     D \in \mathcal{V}$.  $q$ may be $\bdd D$ or it may be another 
     component of $\bdd Q$ but in either case $d(\mathcal{V}, q) \leq 1$.

    Suppose $Q$ contains no disks at all and thus $\chi(Q)\leq 0$.
    Let $E$ be a compressing disk for $N$ in $M$ so that, $|E \cap Q|$
    is minimal among all such disks.  Circles of intersection between
    $Q$ and $E$ and arcs of intersection that are inessential in $Q$
    can be removed by isotoping $E$ via standard innermost disk and
    outermost arc arguments, so this choice of E guarantees that $E$
    and $Q$ only intersect along arcs that are essential in $Q$.  If
    in fact they don't intersect at all, then $d(\bdd E, q) \leq 1$
    for every $q \in \bdd Q$ and we are done.  Consider, then, an arc
    $\bbb$ of $Q \cap E$ that is outermost in $E$, cutting off from
    $E$ a $\bdd$-compressing disk $E_{0}$ for $Q$ that is incident to
    $N$.  Boundary compressing $Q$ along $E_{0}$ gives a new essential
    (by Lemma \ref{lemma:boundredessential}) surface $Q' \subset M$
    which can be isotoped so that each component of $\bdd Q'$ is
    disjoint from each component of $\bdd Q$.  That is for each
    component $q$ of $\bdd Q$ and each component $q'$ of $\bdd Q'$ we
    have that $d(q, q') \leq 1$.
   
     The proof now is by induction on $1-\chi(Q)$.  As $Q$ has no disk
     components, $1-\chi(Q) \geq 1$.  Suppose $1-\chi(Q) = 1$, i.e.
     all components of $Q$ are annuli, so $Q$ is strictly essential.
     As we are not making any claims about the curves in $\mathcal
     {Q}$ coming from $\bdd$-parallel annuli components, we may assume
     all annuli in $Q$ are essential.  Then $Q'$ contains a
     compressing disk $D$ for $N$ (the result of boundary reducing an
     essential annulus component of $Q$ along $E_{0}$) and $\bdd D$ is
     disjoint from all $q \in \bdd Q$.  As $\bdd D \in \mathcal{V}$,
     $d(q, \mathcal{V})\leq 1= 1-\chi(Q)$ as desired.

     Now suppose $1-\chi(Q) > 1$.  If $Q$ is not strictly essential,
     then it contains at least two non-annulus components and, since
     it is essential, at least one essential component.  Thus 
     there is a component $Q_0$ of $Q$ which is essential and
     such that $1-\chi(Q_0) < 1-\chi(Q)$.  By the induction hypothesis,
     for each component $q_0$ of $\bdd Q_0$, $d (q_0, \mathcal{V})\leq
     1-\chi(Q_0)$.  Of course also $d(q,q_0) \leq 1$.  Combining these
     inequalities, we obtain the desired result.

    Suppose next that $Q$ is strictly essential and again all
    $\bdd$-parallel annuli have been removed prior to the boundary
    compression described above.  If the boundary compression creates
    a disk component of $Q'$ then it must be essential and incident to
    $N$ so $\bdd D \in \mathcal{V}$ and for every $q \in \bdd Q$, $d
    (q, \mathcal{V})\leq d(q, \bdd D) \leq 1 \leq 1-\chi(Q)$ and we
    are done.  Suppose then that no component of $Q'$ is a disk and
    $q_{1}$ is any boundary component of an essential component
    $Q_{1}$ of $Q'$.  As $1-\chi (Q_1) \leq 1-\chi (Q') < 1-\chi (Q)$,
    the induction hypothesis applies and $d(q_1,\mathcal{V}) \leq
    1-\chi (Q_1) < 1-\chi (Q)$.  Since for every component $q$ of
    $\bdd Q$, $d(q, q_1) \leq 1$, we have the inequality $d(q,
    \mathcal{V}) \leq d(q_1,\mathcal{V}) + d(q, q_1) \leq
    1-\chi(Q')+1=1-\chi(Q)$, as desired.  \end{proof}

In order to prove Hartshorn's theorem on Heegaard splittings it will
be helpful to understand what it takes to be an essential surface in a
compression body.  Recall the definitions (cf \cite{Sc}):

A {\em compression body} $H$ is a connected 3-manifold obtained from a
closed surface $\bdd_- H$ by attaching 1-handles to $\bdd_- H \times
\{ 1 \} \subset \bdd_- H \times I$.  (It is conventional to consider a
handlebody to be a compression body in which $\bdd_- H = \emptyset$.)
Dually, H is obtained from a connected surface
$\bdd_+ H$ by attaching 2-handles to $\bdd_+ H \times \{ 1 \} \subset
\bdd_+ H \times I$ and 3-handles to any 2-spheres thereby created.
The cores of the $2$-handles are called {\em meridian disks} and a
collection of meridian disks is called {\em complete} if its
complement is $\bdd_- H \times I$, together perhaps with some 
$3$-balls.

Suppose two compression bodies $H_1$ and $H_2$ have $\bdd_+ H_1 \simeq
\bdd_+ H_2$.  Then glue $H_1$ and $H_2$ together along $\bdd_+ H_i =
S$.  The resulting compact 3-manifold $M$ can be written $M = H_1
\cup_{S} H_2$ and this structure is called a {\em Heegaard
splitting} of the 3-manifold with boundary $M$ (or, more specifically,
of the triple ($M; \bdd_- H_1, \bdd_- H_2$)).  It is easy to show that
every compact 3-manifold has a Heegaard splitting.

The following is probably well-known:

\begin{lemma} \label{lemma:compess}
    Suppose $H$ is a compression body and $(Q, \bdd Q) \subset (H,
    \bdd H)$ is incompressible.  If $\bdd Q \cap 
    \bdd_{+}H = \emptyset$, $Q$ is inessential.  That is, each 
    component is $\bdd$-parallel.
    \end{lemma}
    
    \begin{proof} It suffices to consider the case in which $Q$ is
    connected.  To begin with, consider the degenerate case in which
    $H = \bdd_{-} H \times I$.  Suppose there is a counterexample; let
    $Q$ be a counterexample that maximizes $\chi(Q)$.
    
    {\bf Case 1:} $H = \bdd_{-} H \times I$ and $Q$ has non-empty boundary.  
    
    $Q$ cannot be a disk since $\bdd_{-} H \times I$ is
    $\bdd$-irreducible, so $\chi(Q) \leq 0$.  By hypothesis, $\bdd Q
    \subset \bdd_{-} H \times \{ 0 \}$.  Choose $\aaa \subset \bdd_{-}
    H \times \{ 0 \}$ to be any curve that cannot be isotoped off of
    $\bdd Q$ and let $A = \aaa \times I$ be the corresponding annulus
    in $\bdd_{-} H \times I$.  Consider $Q \cap A$ and minimize by
    isotopy of $A$ the number of its components.  A standard argument
    shows that there are no inessential circles of intersection and
    each arc of intersection is essential in $Q$.  Since $\bdd Q$ is
    disjoint from $\bdd_{-} H \times \{ 1 \}$, all arcs of $Q \cap A$
    have both ends in $\bdd_{-} H \times \{ 0 \}$.  An outermost such
    arc in $A$ defines a $\bdd$-compression of $Q$.  The resulting
    surface $Q'$ is still incompressible (for a compressing disk for
    $Q'$ would persist into $Q$) and has at most two components, each
    of higher Euler characteristic and so each $\bdd$-parallel into
    $\bdd_{-} H$.  If there are two components, neither is a disk,
    else the arc along which $\bdd$-compression was supposedly
    performed would not have been essential.  If there are two
    components of $Q'$ and they are not nested (that is, each is
    parallel to the boundary in the complement of the other) it
    follows that $Q$ was $\bdd$-parallel.  If $Q'$ had two nested
    components, it would follow that $Q$ was compressible, a
    contradiction.  (See the end of the proof of Lemma
    \ref{lemma:boundredessential} or Figure \ref{fig:essential}.)
    Similarly, if $Q'$ is connected then, depending on whether the
    tunneling arc dual to the $\bdd$-compression lies inside or
    outside the region of parallelism between $Q'$ and $\bdd M$, $Q$
    would either be compressible or itself $\bdd$-parallel.

    {\bf Case 2:} $H = \bdd_{-} H \times I$ and $Q$ is closed.  
    
    Let $A = \aaa \times I \subset \bdd_{-} H \times I$ be any
    incompressible spanning annulus.  A simple homology argument
    shows that $Q$ intersects $A$.  After the standard move
    eliminating innermost disks, all intersection components are then
    essential curves in $A$.  Let $\lll$ be the curve that is closest
    in $A$ to $\bdd_{-} H \times \{ 0 \}$.  Let $Q'$ be the properly
    embedded surface (now with boundary) obtained from $Q$ by removing
    a neighborhood of $\lll$ in $Q$ and attaching two copies of the
    subannulus of $A$ between $\aaa \times \{ 0 \}$ and $\lll$.  It's
    easy to see that $Q'$ is still incompressible and its boundary is
    still disjoint from $\bdd_{-} H \times \{ 1 \}$, and now $Q'$ has
    non-empty boundary, so by Case 1, $Q'$ is $\bdd$-parallel.  The
    subsurface of $\bdd M$ to which $Q'$ is $\bdd$-parallel can't
    contain the neighborhood $\eta$ of $\aaa \times \{ 0 \}$ in $\bdd
    M$, else the parallelism would identify a compressing disk for
    $Q$.  It follows that the parallelism is outside of $\eta$ and so
    can be extended across $\eta$ to give a parallelism between $Q$
    and a subsurface (hence all) of $\bdd_{-} H \times \{ 0 \}$.
    
    {\bf Case 3:} The general case.  
    
    Let $\Ddd$ be a complete family of meridian disks for $H$, so 
    when $H$ is compressed along $\Ddd$ it becomes a product 
    $\bdd_{-} H \times I$.  Since $Q$ is incompressible, a standard 
    innermost disk argument allows $\Ddd$ to be redefined so that 
    $\Ddd \cap Q$ has no simple closed curves of intersection.  
    Since $Q \cap \bdd_{+} H = \emptyset$ it then follows that $Q 
    \cap \Ddd = \emptyset$.  Then in fact $Q \subset \bdd_{-} H \times 
    I$ and the result follows from Cases 1) or 2).
    \end{proof}

\section{Hartshorn's theorem}

Here we give a quick proof of Hartshorn's theorem (actually, an
extension to the case in which $M$ is not closed) using Proposition
\ref{prop:essential}.  Recall that the distance $d(P)$ of a Heegaard
splitting (\cite{He}) is the minimum distance in the curve complex of
$P$ between a vertex representing a meridian curve on one side of $P$
and a vertex representing a meridian curve on the other side.  

\begin{thm} \label{thm:Hartshorn} Suppose $P$ is a Heegaard splitting
surface for a compact orientable manifold $M$ and $(Q, \bdd Q) \subset (M,
\bdd M)$ is a connected essential surface.  Then $d(P) \leq 2 -
\chi(Q)$.
\end{thm}

Remark: As long as $Q$ contains no inessential disks or spheres, and
at most one essential disk or sphere, $Q$ need not be connected.

\begin{proof} The following are classical facts about Heegaard
splittings (cf \cite{Sc}): If $Q$ is a sphere then $P$ is reducible,
hence $d(P) = 0$.  If $Q$ is a disk then $P$ is $\bdd$-reducible so
$d(P) \leq 1$.  If neither occurs, then $M$ is irreducible and
$\bdd$-irreducible, which is what we henceforth assume.  Moreover,
once $Q$ is neither a disk nor a sphere then $2 - \chi(Q) \geq 2$ so
we may as well assume that $d(P) \geq 2$, ie $P$ is strongly irreducible.

Let $A$, $B$ be the compression-bodies into which $P$ divides $M$ and 
let $\Sss^{A}, \Sss^{B}$ be spines of $A$ and $B$ respectively.  That 
is, $\Sss^{A}$ is the union of a graph in $A$ with $\bdd_{-} A$ and 
$\Sss^{B}$ is the union of a graph in $B$ with $\bdd_{-} B$ so that 
$M - (\Sss^{A} \cup \Sss^{B})$ is homeomorphic to $P \times (-1, 1)$.  
We consider the curves $P \cap Q$ as $P$ sweeps from a neighborhood 
of $\Sss^{A}$ (i. e. near $P \times \{ -1 \}$) to a neighborhood of 
$\Sss^{B}$ (near $P \times \{ 1 \}$).  Under this parameterization, 
let $P_{t}$ denote $P \times \{ t \}$.  Consider the possibilities:

Suppose $Q \cap \Sss^{A} = \emptyset$.  Then $Q$ is an incompressible
surface in the compression body $closure(Q - \Sss^{A}) \cong B$.  By
Lemma \ref{lemma:compess}, $Q$ would be inessential, so this case does
not arise.  Similarly we conclude that $Q$ must intersect $\Sss^{B}$.
It follows that when $t$ is near $-1$, $P_{t} \cap Q$ contains
meridian circles for $A$; when $t$ is near $1$, it contains meridian
circles for $B$.  Since $P$ is strongly irreducible, it can never be
the case that both occur, so at some generic level neither occurs.
(See \cite{Sc} for details, including why we can take such a level to
be generic.)  Hence there is a generic $t_{0}$ so that 
$P_{t_{0}} \cap Q$ contains no meridian circles for $P$.  

An innermost inessential circle of intersection in $P_{t_{0}}$ must be
inessential in $Q$ since $Q$ is incompressible.  So all such circles
of intersection can be removed by an isotopy of $Q$.  After this
process, all remaining curves of intersection are essential in
$P_{t_{0}}$.  Since $P_{t_{0}} \cap Q$ contains no meridian circles
for $P$, no remaining circle of intersection can be inessential in $Q$
either.  Hence all components of $P_{t_{0}} \cap Q$ are essential in
both surfaces; in particular no component of $Q - P_{t_{0}}$ is a
disk.  At this point, revert to $P$ as notation for $P_{t_{0}}$.

If $P \cap Q = \emptyset$ then we are done, just as in the case in
which $Q$ is disjoint from a spine.  Similarly we are done if the
surface $Q_{A} = Q \cap A$ is inessential (hence $\bdd$-parallel) in
$A$ or $Q_{B} = Q \cap B$ is inessential in $B$.  We conclude that 
$Q_{A}$ and $Q_{B}$ are both essential in $A$ and $B$ respectively,
and the positioning of $P$ has guaranteed that no component of either 
is a disk.  

Unless $Q_{A}$ and $Q_{B}$ are both strictly essential, the proof
follows easily from Proposition \ref{prop:essential}: Suppose, for
example, that $Q_{A}$ is not strictly essential and let $\mathcal{U},
\mathcal{V}$ be the set of curves in $P$ bounding disks in $A$ and $B$
respectively.  Let $q$ be a curve in $P \cap Q$ lying on the boundary
of an essential component of $Q_{B}$.  Then Proposition
\ref{prop:essential} says that $d(q, \mathcal{U}) \leq 1 -
\chi(Q_{A})$ and $d(q, \mathcal{V}) \leq 1 - \chi(Q_{B})$ so $$d(P) =
d(\mathcal{U}, \mathcal{V}) \leq d(q, \mathcal{U}) + d(q, \mathcal{V})
\leq (1 - \chi(Q_{A})) + (1 - \chi(Q_{B})) =
2 - \chi(Q)$$ as required.

The case in which $Q_{A}, Q_{B}$ are strictly essential is only a bit
more difficult: Imagine coloring each component of $Q_{A}$ (resp
$Q_{B}$) that is not a $\bdd$-parallel annulus red (resp blue).  Since
$Q_{A}$ and $Q_{B}$ are both essential, there are red and blue regions
in $Q - P$.  Since $Q$ is connected there is a path in $Q$ (possibly
of length $0$) with one end at a red region, one end at a blue region
and no interior point in a colored region.  Since the interior of the
entire path lies in a collection of $\bdd$-parallel annuli, it follows
that the curves in $P \cap Q$ to which the ends of the path are
incident are isotopic curves in $P$.  Now apply the previous argument
to a curve $q \subset P$ in that isotopy class of curves in $P$.
\end{proof}

\bigskip

\section{Sobering examples of large distance} \label{sect:exam}

It is natural to ask whether Proposition \ref{prop:essential} can, in
any useful way, be extended to surfaces that are not essential.  It
appears to be unlikely.  If one allows $Q$ to be $\bdd$-parallel,
obvious counterexamples are easy: take a simple closed curve $\ggg$ in
$N$ that is arbitrarily distant from $\mathcal{V}$ and use for
$Q$ a $\bdd$-parallel annulus $A$ constructed by pushing a regular
neighborhood of $\ggg$ slightly into $M$.  Even if one rules out
$\bdd$-parallel surfaces but does allow $Q$ to be compressible, a
counterexample is obtained by tubing, say, a possibly knotted torus in
$M$ to an annulus $A$ as just constructed.

On the other hand, it has been a recent theme in the study of embedded
surfaces in $3$-manifolds that, for many purposes, a connected
separating surface $Q$ in $M$ will behave much like an incompressible
surface if $Q$ compresses to both sides, but not via disjoint disks. 
Would such a condition on $Q$ be sufficient to guarantee the
conclusion of Proposition \ref{prop:essential}?  That is:
    \begin{quest} \label{quest:main} Suppose $M$ is an irreducible
    compact orientable $3$-manifold, and $N$ is a compressible boundary
    component of $M$. Let $\mathcal{V}$ be the set of 
    essential curves in $N$
    that bound disks in $N$. 
    Suppose further that $(Q, \bdd Q) \subset (M, \bdd M)$ is a
    connected separating surface and $q$ is any boundary component of 
    $Q$.    If $Q$ is compressible into both
    complementary components, but not via disjoint disks, must it be
    true that $d(q, \mathcal{V}) \leq 1-\chi(Q)$?
    
    \end{quest}

In this section we show that there is an example for which the
answer to Question \ref{quest:main} is no.  More remarkably, the next section
shows that it is the only type of bad example.

A bit of terminology will be useful.  
Regard $\bdd D^{2}$ as the end-point union of two arcs, $\bdd_{\pm} D^{2}$.
\begin{itemize}

    \item Suppose $Q \subset M$ is a properly embedded surface and $\ggg \subset
interior(M)$ is an embedded arc which is incident to $Q$
precisely at $\bdd \ggg$.  
There is a relative tubular neighborhood $\eta(\ggg) \cong \ggg \times
D^{2}$ so that $\eta(\ggg)$ intersects $Q$ precisely in the two
disk fibers at the ends of $\ggg$.  Then the surface obtained from $Q$
by removing these two disks and attaching the cylinder $\ggg \times
\bdd D^{2}$ is said to be obtained by {\em tubing along 
\ggg}. 

\item Similarly, suppose $\ggg \subset \bdd M$ is an embedded arc
which is incident to $\bdd Q$ precisely in $\bdd \ggg$. 
There is a relative tubular neighborhood $\eta(\ggg) \cong \ggg
\times D^{2}$ so that $\eta(\ggg)$ intersects $Q$ precisely in the
two $D^{2}$ fibers at the ends of $\ggg$ and $\eta(\ggg)$ intersects
$\bdd M$ precisely in the rectangle $\ggg \times \bdd_{-}D^{2}$. 
Then the properly embedded surface obtained from $Q$ by removing the
two $D^{2}$ fibers at the ends of $\ggg$ and attaching the rectangle $\ggg
\times \bdd_{+}D^{2} $ is said to be obtained by
{\em tunneling along \ggg}.

\end{itemize}

Let $P_{0}, P_{1}$ be two connected compact subsurfaces in the same 
component $N$ of $\bdd M$,
with each component of $\bdd P_{i}, i = 0, 1$ essential in $\bdd M$
and $P_{0} \subset interior(P_{1})$.  Let $Q_{1}$ be the properly
embedded surface in $M$ obtained by pushing $P_{1}$, rel $\bdd$, into
the interior of $M$.  Let $Q_{0}$ denote the properly embedded surface
obtained by pushing $P_{0}$ rel $\bdd$ into the collar between $P_{1}$
and $Q_{1}$.  Then the region $R$ lying between $Q_{0}$ and $Q_{1}$ is
naturally homeomorphic to $Q_{1} \times I$.  (Here $\bdd Q_{1} \times
I$ can be thought of either as vertically crushed to $\bdd Q_{1}
\subset \bdd M$ or as constituting a small collar of $\bdd Q_{1}$ in
$P_{1} \subset \bdd M$.)  Under the homeomorphism $R \cong Q_{1}
\times I$ the top of $R$ (corresponding to $Q_{1} \times \{1 \}$) is
$Q_{1}$ and the bottom of $R$ (corresponding to $Q_{1} \times \{0 \}$)
is the boundary union of $Q_{0}$ and $P_{1} - P_{0}$.  The properly
embedded surface $Q_{0} \cup Q_{1} \subset M$ is called the {\em
recessed collar} determined by $P_{0} \subset P_{1}$ bounding $R$.

Recessed collars behave predictably under tunnelings:

\begin{lemma} \label{lemma:reccollar}
    
    Suppose $Q_{0} \cup Q_{1} \subset M$ is the recessed collar
    determined by $P_{0} \subset P_{1}$, and $R \cong P_{1} \times I$
    is the component of $M - Q$ on whose boundary both $Q_{i}$ lie.
    Let $\ggg \subset \bdd M$ be a properly embedded arc in $\bdd M -
    (Q_{0} \cup Q_{1})$.  Let $Q_{+}$ be the surface obtained from
    $Q_{0} \cup Q_{1}$ by tunneling along $\ggg$.  Then
    
    \begin{enumerate}
	
	\item If $\ggg \subset (P_{1} - P_{0})$ and $\ggg$ has both
	ends on $\bdd P_{0}$ or if $\ggg \subset (\bdd M - P_{1})$,
	then $Q_{+}$ is a recessed collar.
	
	\item If $\ggg \subset P_{0}$ then there is a compressing disk
	for $Q_{+}$ in $M - R$.
	
	\item If $\ggg \subset (P_{1} -
	P_{0})$ and $\ggg$ has one or both ends on $\bdd P_{1}$, then
	there is a compressing disk for $Q_{+}$ in $R$.
	
	\end{enumerate}
	
\end{lemma}
	
\begin{proof} In the first case, tunneling is equivalent to just
adding a band to either $P_{1}$ or $P_{0}$ and then constructing the
recessed collar.  In the second case, the disk $\ggg \times I$ in the
collar between $P_{0}$ and $Q_{0}$ determines a compressing disk for
$Q_{+}$ (that is, for the component of $Q_{+}$ coming from $Q_{0}$)
that lies outside of $R$.  

Similarly, in one of the third cases, when $\ggg \subset (P_{1} -
P_{0})$ has both ends on $\bdd P_{1}$, $\ggg \times I$ in the collar
between $P_{1}$ and $Q_{1}$ determines a compressing disk for $Q_{+}$
(this time for the component of $Q_{+}$ coming from $Q_{1}$) that this
time lies inside of $R$.  

In the last case, when one end of $\ggg \subset (P_{1} - P_{0})$ lies
on each of $\bdd P_{0}$ and $\bdd P_{1}$ a slightly more sophisticated
construction is needed.  After the tunneling construction, $\bdd Q_{+}
\cap interior(P_{1})$ has one arc component $\ggg'$ that consists of two
parallel copies of the spanning arc $\ggg$ and a subarc of the 
component of $\bdd P_{0}$ that is incident to $\ggg$.  This arc 
$\ggg' \subset \bdd Q_{+}$ can be pushed slightly into $Q_{+}$.  Then 
the disk $\ggg' \times I$ (using the product structure on $R$) 
determines a compressing disk for $Q_{+}$ that lies in $R$.  
(The disk $\ggg' \times I$ looks much like the disk $E$ in Figure 
\ref{fig:essential}.)
\end{proof}

One of the constructions of this lemma will be needed in a different context:

\begin{lemma}  \label{lemma:2collars}
    Suppose $Q_{0} \cup Q_{1} \subset M$ and $Q_{1} \cup Q_{2} \subset
    M$ are the recessed collars determined by connected surfaces $P_{0}
    \subset interior(P_{1})$ and $P_{1} \subset interior(P_{2})$.  Let
    $R_{1}, R_{2}$ be the regions these recessed collars bound. 
    Furthermore, let $\ggg_{i} \subset \bdd M, i = 1, 2$ be properly
    embedded arcs spanning $P_{1} - P_{0}$ and $P_{2} - P_{1}$
    respectively.  That is, $\ggg_{i}$ has one end point on each of
    $\bdd P_{i}, \bdd P_{i-1}$.  Let $Q_{+}$ be the connected surface
    obtained from $Q_{0} \cup Q_{1} \cup Q_{2}$ by tunneling along
    both $\ggg_{1}$ and $\ggg_{2}$.  Then either
    
    \begin{enumerate}
	
	\item There are disjoint compressing disks for $Q_{+}$ in $R_{1}$ and 
	$R_{2}$ or
	
	\item $P_{0}$ is an annulus parallel in $P_{1}$ to a component
	$c$ of $\bdd P_{1}$, and $c$ is incident to both tunnels.
	
	\end{enumerate}
	
	In the latter case, $Q_{+}$ is properly isotopic to the surface
	obtained from the recessed collar $Q_{1} \cup Q_{2}$ by tubing along 
	an arc in $interior(M)$ that is parallel to $\ggg_{2} \subset \bdd M$.
	
\end{lemma}
    
\begin{figure}[tbh]
\centering
\includegraphics[width=0.8\textwidth]{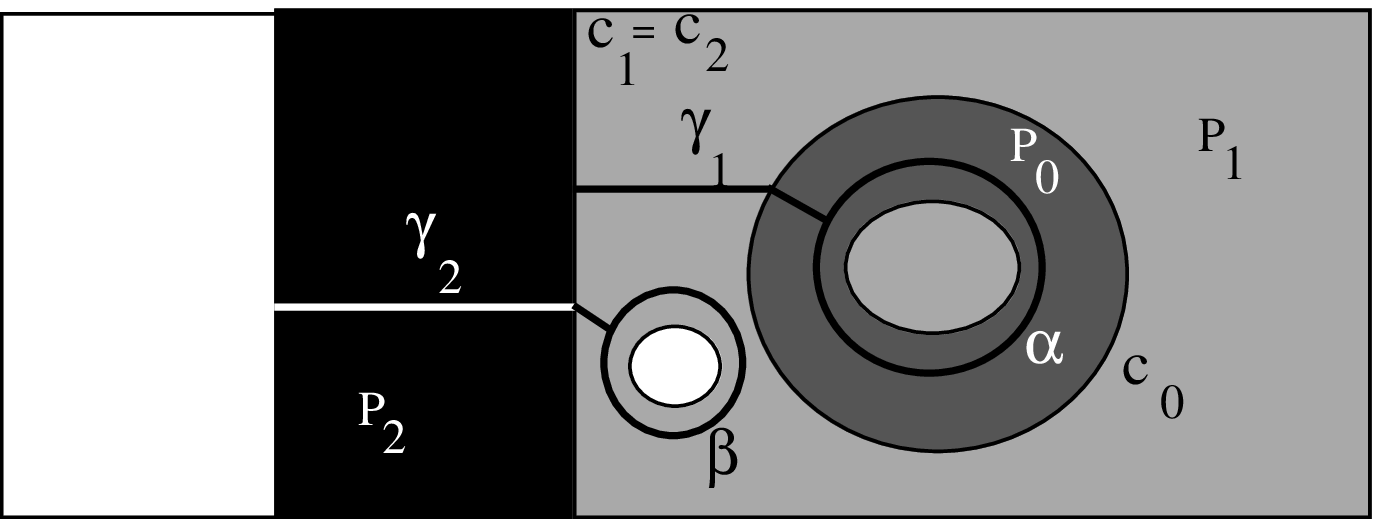}
\caption{} \label{fig:twocollars}
\end{figure}

\begin{proof} For $P$ any surface with boundary, define an {\em
eyeglass graph} in $P$ to be the union of an essential simple closed
curve in the interior of $P$ and an embedded arc in the curve's
complement, connecting the curve to $\bdd P$.  

Let $c_{1} \subset \bdd P_{1}$ and $c_{0} \subset \bdd P_{0}$ be the
components to which the ends of $\ggg_{1}$ are incident.  Let $c_{2}$
be the component of $\bdd P_{1}$ (note: not $\bdd P_{2}$) to which the
end of $\ggg_{2}$ is incident.  (Possibly $c_{1} = c_{2}$.)  Let
$\aaa$ be any essential simple closed curve in $P_{0}$ and choose an
embedded arc in $P_{0} - \aaa$ connecting $\aaa$ to the end of
$\ggg_{1}$ in $c_{0}$; the union of that arc, the closed curve $\aaa$
and the arc $\ggg_{1}$ is an eyeglass curve $e_{1}$ in $P_{1}$ which
intersects $P_{1} - P_{0}$ in the arc $\ggg_{1}$.  Then the
construction of Lemma \ref{lemma:reccollar}, there applied to the
eyeglass $\ggg_{1} \cup c_{0}$, shows here that a neighborhood of the
product $e_{1} \times I \subset R_{1} \cong P_{1} \times I$ contains a
compressing disk for $Q_{+}$ that lies in $R_{1}$ and which intersects
$Q_{1}$ in a neighborhood of $e_{1} \times \{ 1 \}$.

Similarly, for $\bbb$ any essential simple closed curve in $P_{1}$,
and an embedded arc in $P_{1} - \bbb$ connecting $\bbb$ to the end of
$\ggg_{2}$ in $c_{2}$ we get an eyeglass $e_{2} \subset P_{2}$ and a compressing
disk for $Q_{+}$ that lies in $R_{2}$ and whose boundary
intersects $Q_{1}$ only within a neighborhood of $e_{2} \times \{ 1 \}$. 
So if we can find disjoint such eyeglasses in $P_{1}$ and $P_{2}$ we will have
constructed the required disjoint compressing disks.

Suppose first that $P_{0}$ is not an annulus parallel to
$c_{1}$.  Then $P_{0}$ contains an essential simple closed curve $\aaa$ that is
not parallel to $c_{1}$.  Since $\aaa$ is not parallel to $c_{1}$, no 
component of the complement $P_{1} - e_{1}$ is a disk, so there is an 
essential simple closed curve $\bbb$ in the component of $P_{1} - e_{1}$ 
that is incident to $c_{2}$.  The same is true even if $P_{0}$ is an 
annulus parallel to $c_{1}$ so long as $c_{1} \neq c_{2}$.  This 
proves the enumerated items.  See Figure \ref{fig:twocollars}

The proof that in case 2), $Q_{+}$ can be described by tubing $Q_{1}$ 
to $Q_{2}$ along an arc parallel to $\ggg_{2}$ is a pleasant exercise 
left to the reader.    \end{proof}

Now consider a particular type of tubing of a recessed collar.
Suppose $Q_{0} \cup Q_{1} \subset M$ is the recessed collar bounding
$R$ determined by $P_{0} \subset P_{1} \subset \bdd M$.  Let $\rho$
denote a vertical spanning arc in $R$, that is, the image in $R \cong
P_{1} \times I$ of $point \times I$, where $point \in P_{0}$.  Let $Q$
be the surface obtained from $Q_{0} \cup Q_{1}$ by tubing along
$\rho$.  Then $Q$ is called a {\em tube-spanned recessed collar}.

A tube-spanned recessed collar has nice properties:

\begin{lemma} \label{lemma:tubespan}
    
    Suppose $Q$ is a tube-spanned recessed collar constructed as
    above.  Then
    
    \begin{itemize}
	
	\item $Q$ is connected and separating and $Q$ compresses in both  
	complementary components in $M$..
	
	\item If $Q$ compresses in both complementary components via 
	disjoint disks, then $P_{1} \subset \bdd M$ is
	compressible in $M$.
	
	\item If $Q_{+}$ is obtained from $Q$ by tunneling, then
	either $Q_{+}$ is also a tube-spanned recessed collar or
	$Q_{+}$ compresses in both complementary components via
	disjoint disks.  (Possibly both are true).
	
	\item If $Q_{+}$ is obtained from $Q$ by tunneling together
	$Q$ and a $\bdd$-parallel connected incompressible surface
	$Q'$, then either $Q_{+}$ is also a tube-spanned recessed
	collar or $Q_{+}$ compresses in both complementary components
	via disjoint disks.  (Possibly both are true).
	
	\end{itemize}
	
\end{lemma} 

\begin{proof}  The construction guarantees that $Q$ is connected and
separating.  It compresses on both sides: Let $Y$ denote the component
$R - \eta(\rho)$ of $M - Q$ and let $X$ be the other component.  A
disk fiber $\mu$ of $\eta(\rho)$ is a compressing disk for $Q$ in $X$. 
To see a compressing disk for $Q$ in $Y$, start with an
essential simple closed curve in $Q_{0}$ containing 
the end of $\rho$ in $Q_{0}$.  The corresponding vertical annulus $A 
\subset R$ includes the vertical arc $\rho \subset R$.  Then $A -
\eta(\rho)$ is a disk in $Y$ whose boundary is essential in $Q$.

To prove the second property, suppose that there are disjoint
compressing disks, $D_{X} \subset X$ and $D_{Y} \subset Y$.  $\bdd
D_{Y}$ cannot be disjoint from the meridian $\mu$ of $\eta(\rho)$
since if it were, $\bdd D_{Y}$ would lie in either on the top or the
bottom of $Y \cong (P_{1} - point) \times I$, either of which is
clearly incompressible in $Y$.  So $D_{X}$ cannot be parallel to
$\mu$.  A standard innermost disk argument allows us to choose $D_{X}$
so that $D_{X} \cap \mu$ contains no circles of intersection, and an
isotopy of $\bdd D_{X}$ on $Q$ ensures that any arc component of $\bdd
D_{X} - \mu$ is essential in one of the punctured surfaces $Q_{1} \cap
Q$ or $Q_{0} \cap Q$.  If $D_{X}$ is disjoint from $\mu$ it lies on
$Q_{1}$, say, but in any case it determines a compressing disk for
$P_{1}$ in $M$, as required.  If $D_{X}$ is not disjoint from $\mu$
then an outermost disk in $D_{X}$ cut off by $\mu$ would similarly
determine a compression of $P_{1}$ in $M$.

The third property follows from Lemma \ref{lemma:reccollar}.  When the
tunneling there leaves $Q_{+}$ as a recessed collar (option 1) then
the operation here leaves $Q_{+}$ a tube-spanned recessed collar.  If
the tunneling arc $\ggg$ lies in $P_{1} - P_{0}$ and thereby gives
rise to a compressing disk in $R$ (option 3), the compressing disk
$D_{Y}$ there constructed lies in $Y$ and so can clearly be kept
disjoint from the vertical arc $\rho$.  Then $D_{Y}$ is disjoint from
the compressing disk $\mu$ for $X$, as required.  Finally, if $\ggg$
lies in $P_{0}$ then the compressing disk $D_{X}$ in $M - R$
constructed there lies in $X$ and intersects $Q_{0}$ in a single 
essential arc.  The simple closed curve in $Q_{0}$ from which $A$ is 
constructed can be taken to intersect $D_{X}$ in at most one point, 
so in the end the disk $D_{Y} \subset Y$ intersects $D_{X}$ in at most 
one point.  Then the boundary of a regular neighborhood of $\bdd X 
\cup \bdd Y$ in $Q$ is a simple closed curve that bounds a disk in 
both $X$ and $Y$, as required.

The fourth property is proven in a similar way.  Suppose first that
$\bdd Q'$ is disjoint from $P_{1}$.  If the region $P' \subset \bdd M$
to which $Q'$ is parallel is disjoint from $P_{1}$ then tunneling $Q'$
to $Q_{1}$ just creates a larger $\bdd$-parallel surface and $Q_{+}$
is a tube-spanned recessed collar.  If $P_{1} \subset P'$ then the
region $R'$ between $Q'$ and $Q_{1}$ is a recessed collar and
according to option 3 of Lemma \ref{lemma:reccollar} there is a
compressing disk for $Q_{+}$ in $R' \cap X$ that is incident to
$Q_{1}$ only in a collar of $\bdd Q_{1}$.  In particular it is
disjoint from a compressing disk for $Q$ in $R \cap Y$ constructed
above from an annulus $A$ that is incident to $Q_{1}$ away from this
collar.

Next suppose that $\bdd Q'$ lies in $P_{1} - P_{0}$, so $P' \subset
P_{1} - P_{0}$.  If the tunnel connects $Q'$ to $Q_{0}$ then tunneling 
$Q_{0}$ to $Q'$ just creates a larger $\bdd$-parallel surface and
$Q_{+}$ is a tube-spanned recessed collar.  If the tunneling connects 
$Q'$ to $Q_{1}$ then the argument is the same as when $Q_{+}$ is
obtained from $Q$ by tunneling into $P_{1} - P_{0}$ with both ends of 
the tunnel on $\bdd P_{1}$.  

Finally suppose that $\bdd Q'$ lies in $P_{0}$, so $P' \subset P_{0}$.
Then the tunneling connects $Q'$ to $Q_{0}$.  The region $R'$ between
$Q'$ and $Q_{0}$ is a recessed collar and according to option 3 of
Lemma \ref{lemma:reccollar} there is a compressing disk for $Q_{+}$ in
$R' \cap X$ that is incident to $Q'$ only in a collar of $\bdd Q'$.
In particular it is disjoint from the compressing disk for $Q$ in $R
\cap Y$ constructed above from an annulus $A$ incident to $Q_{0}$ in
the image of $P' \subset P_{1}$ away from that collar.
\end{proof}

\begin{cor} 
Suppose $M$ is an irreducible compact orientable $3$-manifold, and $N$
is a compressible boundary component of $M$.  Let $\mathcal{V}$ be the
set of curves in $N$ that arise as boundaries of compressing disks of
$N$.  Then for any $n \in \mathbb{N}$ there is a connected properly
imbedded separating surface $(Q, \bdd Q) \subset (M, N)$ so that $Q$
compresses in both complementary components but not via disjoint
disks and, for any component $q$ of $\bdd Q$, $d(q, \mathcal{V}) \geq n$.
   \end{cor}
   
   \begin{proof}  Let $A_{1}$ be an annulus in $\bdd M$ whose core has 
   distance at least $n$ from $\mathcal{V}$.  Let $A_{0} \subset A_{1}$ be a 
   thinner subannulus and let $Q$ be the tube-spanned recessed 
   product in $M$ that they determine.  The result follows from the first two 
   conclusions of Lemma \ref{lemma:tubespan}.
   
   \end{proof}

\section{Any example is a tube-spanned recessed collar}

It will be useful to expand the context beyond connected 
separating surfaces.    

\begin{defin} \label{defin:splitting}
    
    Let $(Q, \bdd Q) \subset (M, \bdd M)$ be a properly embedded
    orientable surface in the orientable irreducible $3$-manifold $M$.
    $Q$ will be called a {\em splitting surface} if no component is
    closed, no component is a disk, and $M$ is the
    union of two $3$-manifolds $X$ and $Y$ along $Q$.  
    
    We abbreviate by saying that $Q$ splits $M$ into the submanifolds
    $X$ and $Y$. \end{defin}
    
    The definition differs slightly from that of \cite[Definition
    1.1]{JS}, which allows $Q$ to have closed components and disk
    components.
    Note also that the condition that $M$ is the union of two
    $3$-manifolds $X$ and $Y$ along $Q$ is equivalent to saying that
    $Q$ can be normally oriented so that any oriented arc in $M$
    transverse to $Q$ alternately crosses $Q$ in the direction
    consistent with the normal orientation and then against the normal
    orientation.
   
\begin{defin} \label{defin:compressiontypes}
    
Suppose as above that $(Q, \bdd Q) \subset (M, \bdd M)$ is a splitting
surface that splits $M$ into submanifolds $X$ and $Y$.  $Q$ is
bicompressible if both $X$ and $Y$ contain compressing disks for $Q$
in $M$; $Q$ is strongly compressible if there are such disks whose
boundaries are disjoint in $Q$.  If $Q$ is not strongly compressible
then it is weakly incompressible.  \end{defin}
	
	Note that if $Q$ is bicompressible but weakly incompressible $\bdd
Q$ is necessarily essential in $\bdd M$, for otherwise an innermost
inessential component would bound a compressing disk for $Q$ in $Y 
\cap \bdd M$ (say).  Such a disk, lying in $\bdd M$, would necessarily 
be disjoint from any compressing disk for $Q$ in $X$. 
	
There are natural extensions of these ideas.  One extension that will
eventually prove useful is to $\bdd$-compressions of splitting
surfaces:

\begin{defin} \label{defin:boundtypes}
    
	A splitting surface $(Q, \bdd Q) \subset (M, \bdd M)$ is
	strongly $\bdd$-compressible if there are $\bdd$-compressing
	disks $D_{X} \subset X, D_{Y} \subset Y$ and $\bdd D_{X} \cap
	\bdd D_{Y} = \emptyset$.

\end{defin}

Here is our main result:

\begin{thm} \label{thm:bicompressible}  
     Suppose $M$ is an irreducible compact orientable $3$-manifold,
     $N$ is a compressible boundary component of $M$ and $(Q, \bdd Q)
     \subset (M, \bdd M)$ is a bicompressible, weakly incompressible
     splitting surface with a bicompressible component incident to $N$.
     
      Let $\mathcal{V}$ be the set of essential curves in $N$
     that bound disks in $M$ and let $q$ be any component of $\bdd Q
     \cap N$.   Then
     either
     \begin{itemize}
     \item $d(q, \mathcal{V}) \leq 1 - \chi(Q)$ in the curve complex
     on $N$ or 
     \item $q$ lies in the boundary of a $\bdd$-parallel
     annulus component of $Q$ or 
     \item one component of $Q$ is a
     tube-spanned recessed collar; all other components incident to
     $N$ are incompressible and $\bdd$-parallel.
     
     \end{itemize}   
\end{thm}

Note that in the last case, $Q$ lies entirely in a collar of $N$.

\begin{lemma}  \label{lemma:Ddisjoint} 
    
    Let $(Q, \bdd Q) \subset (M, \bdd M)$ be as in Theorem 
    \ref{thm:bicompressible}, splitting $M$ into $X$ and
    $Y$.  Let $Q_{X}$ be the result of maximally compressing $Q$ into
    $X$.  Then
    
    \begin{enumerate}
	\item $Q_{X}$ is incompressible in $M$ and,
	
	\item there is a compressing disk $D$ for $N$ in $M$, so that
	some complete set of compressing disks for $Q$ in $X$ is
	disjoint from $D$ and, moreover, $Q \cap D$ consists entirely
	of arcs that are essential in $Q_{X}$.
	
	\end{enumerate}

\end{lemma}

\begin{proof} First we show that $Q_{X}$ is incompressible.  This is
in some sense a classical result, going back to Haken.  A more modern
view is in \cite{CG}.  Here we take the viewpoint first used in
\cite[Prop.  2.2]{ST}, which adapts well to other contexts we will
need as well and is a good source for details missing here.

$Q_{X}$
is obtained from $Q$ by compressing into $X$.  Dually, we can think of
$Q_{X}$ as a surface splitting $M$ into $X'$ and $Y'$ (except possibly
$Q_{X}$ has some closed components) and $Q$ is constructed from
$Q_{X}$ by tubing along a collection of arcs in $Y'$.  Sliding one of
these arcs over another or along $Q_{X}$ merely moves $Q$ by an
isotopy, so an alternate view of the construction is this: There is a
graph $\Ggg \subset Y'$, with all of its valence-one vertices on $Q_{X}$.  A
regular neighborhood of $Q_{X} \cup \Ggg$ has boundary consisting of a
copy of $Q_{X}$ and a copy of $Q$.  (This construction of $Q$ from
$Q_{X}$ could be called $1$-surgery along the graph $\Ggg$.)  The
graph $\Ggg$ may be varied by slides of edges along other edges or along
$Q_{X}$; the effect on $Q$ is merely to isotope it in the complement 
of $Q_{X}$.

Suppose that $F$ is a compressing disk for $Q_{X}$ in $M$.  $F$ must 
lie in $Y'$, else $Q$ could be further compressed into $X$.  
Choose a representation of \Ggg\ which minimizes $|F \cap \Ggg|$,
and then choose a compressing disk $E$ for $Q$ in $Y$ which
minimizes $|F \cap E|$.  If there are any closed components of $F \cap
E$, an innermost one in $E$ bounds a subdisk of $E$ disjoint from $F$,
$\Ggg$ and $Q$; an isotopy of $F$ will remove the intersection curve
without raising $|F \cap \Ggg|$.  So in 
fact there are no closed curves in $F \cap E$.  

The disk $F$ must intersect the graph $\Ggg$ else $F$ would lie
entirely in $Y$ and so be a compressing disk for $Q$ in $Y$ that is
disjoint from compressing disks of $Q$ in $X$.  This would contradict
the weak incompressiblity of $Q$.  One can view the intersection of
$\Ggg \cup E$ with $F$ as a graph $\Lll \subset F$ whose vertices are
the points $\Ggg \cap F$ and whose edges are the arcs $F \cap E$.

If there is an isolated vertex of the graph $\Lll \subset F$ (i.  e. a
point in $\Ggg \cap F$ that is disjoint from $E$) then the vertex
would correspond to a compressing disk for $Q$ in $X$ which is
disjoint from $E$, contradicting weak irreducibility.  If there is a
loop in $\Lll \subset F$ whose interior contains no vertex, an
innermost such loop would bound a subdisk of $F$ that could be used to
simplify $E$; that is to find disk $E_{0}$ for $Q$ in $Y$ so that $|F
\cap E_{0}| < |F \cap E|$, again a contradiction.  We conclude that
$\Lll$ has a vertex $w$ that is incident to edges but to no loops of
$\Lll$.  Choose an arc \bbb\ which is outermost in $E$ among all arcs
of $F \cap E$ which are incident to $w$.  Then \bbb\ cuts off from $E$
a disk $E'$ with $E'-\bbb$ disjoint from $w$.  Let $e$ be the edge of
\Ggg\ which contains $w$.  Then the disk $E'$ gives instructions about
how to isotope and slide the edge $e$ until $w$ and possibly other
points of $\Ggg \cap F$ are removed, lowering $|\Ggg \cap F|$, a
contradiction that establishes the first claim.

To establish the second claim, first note that by shrinking a complete
set of compressing disks for $Q$ in $X$ very small, we can of course
make them disjoint from any $D$; the difficulty is ensuring that $Q_{X}
\cap D$ then has no simple closed curves of intersection.

Choose $D$ and isotope $Q_{X}$ to minimize the number of components
$|D \cap Q_{X}|$, then choose a representation of \Ggg\ which
minimizes $|D \cap \Ggg|$, and finally then choose a compressing disk
$E$ for $Q$ in $Y$ which minimizes $|D \cap E|$.  If there are any
closed components of $D \cap E$, an innermost one in $E$ bounds a
subdisk of $E$ disjoint from $D, \Ggg$ and $Q$; an isotopy of $D$ will
remove the intersection curve without raising either $|D \cap Q_{X}|$
or $|D \cap \Ggg|$.  So in fact there are no closed curves in $D \cap
E$.

Suppose there are closed curves in $D \cap Q_{X}$.  An innermost one
in $D$ will bound a subdisk $D_{0}$.  Since $Q_{X}$ is incompressible,
$\bdd D_{0}$ also bounds a disk in $Q_{X}$; the curve of intersection
could then be removed by an isotopy of $Q_{X}$, a 
contradiction.

From this contradiction we deduce that all components of $D \cap
Q_{X}$ are arcs.  All arcs are essential in $Q_{X}$ else $|D \cap
Q_{X}|$ could be lowered by rechoosing $D$.  The only other components
of $D \cap Q$ are closed curves, compressible in $X$, each
corresponding to a point in $D \cap \Ggg$.  So it suffices to show
that $D \cap \Ggg = \emptyset$.  The proof is analogous to the proof
of the first claim, where it was shown that $\Ggg$ must be
disjoint from any compressing disk $F$ for $Q_{X}$ in $Y'$, but now
for $F$ we take a (disk) component of $D - Q_{X}$.

If no component of $D - Q_{X}$ intersects $\Ggg$ there is nothing to
prove, so let $F$ be a component intersecting $\Ggg$ and regard $\Lll
= (\Ggg \cup E) \cap F$ as a graph in $F$, with possibly some edges
incident to the arcs $Q_{X} \cap D$ lying in $\bdd F$.  As above, no
vertex of $\Lll$ (i.  e. point of $\Ggg \cap F$) can be isolated in
$\Lll$ and an innermost inessential loop in $\Lll$ would allow an
improvement in $E$ so as to reduce $D \cap E$.  Hence there is a
vertex $w$ of $\Lll$ that is incident to edges but no loops in $\Lll$.
An edge in $\Lll$ that is outermost in $E$ among all edges incident to
$w$ will cut off a disk from $E$ that provides instructions how to
slide the edge $e$ of $\Ggg$ containing $w$ so as to remove the
intersection point $w$ and possibly other intersection points.  As in
the first case, some sliding of the end of $e$ may necessarily be
along arcs in $Q_{X}$, as well as over other edges in $\Ggg$.  
	\end{proof}
	
	\bigskip
	
{\bf Proof of Theorem \ref{thm:bicompressible}:}  Just as in
the proof of Proposition \ref{prop:essential} the proof is by
induction on $1 - \chi(Q)$.  Since $Q$ contains no disk components, 
$1 - \chi(Q) \geq 1$.  

If compressing disks for $Q$ were incident to two different components
of $Q$, then there would be compressing disks on opposite sides
incident to two different components of $Q$, violating weak
incompressibility.  So we deduce that all compressing disks for $Q$
are incident to at most one component $Q_{0}$ of $Q$.  $Q_{0}$ cannot
be an annulus, else the boundaries of compressing disks in $X$ and $Y$
would be parallel in $Q_{0}$ and so could be made disjoint.  If $Q$
also contains an essential component $Q'$ incident to $N$ then $1 -
\chi(Q') \leq 1 - \chi(Q - Q_{0}) < 1 - \chi(Q)$ and so, by
Proposition \ref{prop:essential}, for any component $q'$ of $\bdd Q'
\cap N$, $d(q', \mathcal{V}) \leq 1 - \chi(Q') < 1 - \chi(Q)$.  This
implies that $d(q, \mathcal{V}) \leq d(q', \mathcal{V}) + d(q, q')
\leq 1 - \chi(Q)$ as required.  So we will also henceforth assume that
no component of $Q$ incident to $N$ is essential.

We can also assume that each component of $Q - Q_{0}$ is itself an
incompressible surface.  For suppose $D$ is a compressing disk for a
component $Q_{1} \neq Q_{0}$ of $Q$, chosen among all such disks to
have a minimal number of intersection components with $Q$.  If the
interior of $D$ were disjoint from $Q$ then $D$ would be a compressing
disk for $Q$ itself, violating weak incompressibility as described
above.  Similarly, an innermost circle of $Q \cap D$ in $D$ must lie
in $Q_{0}$.  Consider a subdisk $D'$ of $D$ (possibly all of $D$) with
the property that its boundary is second-innermost among components of
$D \cap Q$.  That is, the interior of $D'$ intersects $Q$ exactly in
innermost circles of intersection, each bounding disks in $X$, say.
If $\bdd D'$ is not in $Q_{0}$ then it is also a compressing disk for
$Q_{X}$, contradicting the first statement in Lemma
\ref{lemma:Ddisjoint}.  The argument is only a bit more subtle when
$\bdd D'$ is in $Q_{0}$, cf the No Nesting Lemma 
\cite[Lemma 2.2]{Sc2}. 

Let $Q_{-}$ be the union of components of $Q$ that are not incident to
$N$.  Since $Q_{-}$ is incompressible, each compressing disk for $N$
is disjoint from $Q_{-}$.  In particular, it suffices to work inside
the $3$-manifold $M - \eta(Q_{-})$ instead of $M$.  So, with no loss
of generality, we can assume that $Q_{-} = \emptyset$, ie each
component of $Q$ is incident to $N$.  

Since each component of $Q$ other than $Q_{0}$ is
incompressible and not essential, each is boundary parallel.  In
particular, removing one of these components $Q_{1}$ from $Q$ still
leaves a bicompressible, weakly incompressible splitting surface,
though each component of $M - Q_{1}$ in the region of parallelism
between $Q_{1}$ and $\bdd M$ would need to be switched from $X$ to $Y$
or vice versa.  Since we don't care about the boundaries of
$\bdd$-parallel annuli, all such components can be removed from $Q$
without affecting the hypotheses or conclusion.  If there remains a
$\bdd$-parallel component $Q_{1}$ that is not an annulus, then
consider $Q' = Q - Q_{1}$.  We have $1 - \chi(Q') < 1 - \chi(Q)$ so
the inductive hypothesis applies.  Then either $Q_{0}$ is a
tube-spanned recessed product (and we are done) or for any component
$q'$ of $\bdd Q'$, $d(q', \mathcal{V}) \leq 1 - \chi(Q') < 1 -
\chi(Q)$.  This implies that $d(q, \mathcal{V}) \leq d(q',
\mathcal{V}) + d(q, q') \leq 1 - \chi(Q)$ and again we are done.  So
we may as well assume that $Q = Q_{0}$ is connected and, as we have
seen, not an annulus.

\bigskip

{\bf Claim:}  The theorem holds if $Q$ is strongly $\bdd$-compressible.

{\bf Proof of claim:} Suppose there are disjoint $\bdd$-compressing
disks $F_{X} \subset X$, $F_{Y} \subset Y$ for $Q$ in $M$.  Let
$Q_{x}, Q_{y}$ denote the surfaces obtained from $Q$ by
$\bdd$-compressing $Q$ along $F_{X}$ and $F_{Y}$ respectively, and let
$Q_{-}$ denote the surface obtained by $\bdd$-compressing along both
disks simultaneously.  (We use lower case $x, y$ to distinguish these
from the surfaces $Q_{X}, Q_{Y}$ obtained by {\em maximally
compressing} $Q$ into respectively $X$ or $Y$.)  A standard innermost
disk, outermost arc argument between $F_{x}$ and a compressing disk
for $Q$ in $X$ shows that $Q_{x}$ is compressible in $X$.  Similarly,
$Q_{y}$ is compressible in $Y$.

Each of $Q_{x}, Q_{y}$ has at most two components, since $Q$ is
connected.  Suppose that $Q_{x}$ (say) is itself bicompressible.  If
it were strongly compressible, the same strong compression pair of
disks would strongly compress $Q$, so we conclude that the inductive
hypothesis applies to $Q_{x}$, so we apply the theorem to $Q_{x}$.
One possibility is that one component of $Q_{x}$ is a tube-spanned
recessed collar and the other (if there are two components) is
$\bdd$-parallel.  But by Lemma \ref{lemma:tubespan} this case implies
that $Q$ is also a tube-spanned recessed collar and we are done.  The
other possibility is that for $q_{x}$ a component of the boundary of
an essential component of $Q_{x}$, $d(q_{x}, \mathcal{V}) \leq
1-\chi(Q_{x}) < 1-\chi(Q)$.  This implies that $d(q, \mathcal{V}) \leq
d(q_{x}, \mathcal{V}) + d(q, q_{x}) \leq 1 - \chi(Q)$ and again we are
done.  So we henceforth assume that $Q_{x}$ (resp $Q_{y}$) is
compressible into $X$ (resp $Y$) but not into $Y$ (resp $X$).

It follows that $Q_{-}$ is incompressible, for if $Q_{-}$ is
compressible into $Y$, say, then such a compressing disk would be
unaffected by the tunneling that recovers $Q_{x}$ from $Q_{-}$ and
$Q_{x}$ would also compress into $Y$.  

On the other hand, if $Q_{-}$ is essential in $M$ then the claim
follows from Proposition \ref{prop:essential}.  So the only remaining
case to consider in the proof of the claim is when $Q_{-}$ is
incompressible and not essential, so all its components are
$\bdd$-parallel.  Since $Q$ is connected, $Q_{-}$ has at most three
components.  Suppose there are exactly three $Q_{0}, Q_{1}, Q_{2}$.
If the three are nested (that is, they can be arranged as $Q_{0},
Q_{1}, Q_{2}$ are in Lemma \ref{lemma:2collars}) then that lemma shows
that the weakly incompressible $Q$ must be a tube-spanned recessed
collar, as required.  If no pairs of the three components of $Q_{-}$
are nested, then $Q$ itself would be boundary parallel and so could
not be compressible on the side towards $N$.  Finally, suppose that
two components ($Q_{0}, Q_{1}$, say) are nested, that $Q_{2}$ is
$\bdd$-parallel in their complement, and $Q_{x}$, say, is obtained
from $Q_{1}, Q_{2}$ by tunneling between $Q_{1}$ and $Q_{2}$, so
$Q_{x}$ is $\bdd$-parallel.  $Q_{x}$ is also compressible; the
compressing disk either also lies in a collar of $N$, or, via the
parallelism to the boundary, the disk represents a compressing disk
$D$ for $N$ in $M$ whose boundary is disjoint from $\bdd Q_{x}$.  In
the latter case we have, for $q_{x}$ any component of $\bdd Q_{x}$,
$d(q_{x}, \bdd D) \leq 1$.  Then for $q$ any component of $Q$, $d(q,
\bdd D) \leq d(q_{x}, \bdd D) + d(q, q_{x}) \leq 2 \leq  1 - \chi(Q)$
and we are done.  The former case can only arise if there are boundary
components of $Q_{1}$ and $Q_{2}$ that cobound an annulus, and that
annulus is spanned by the tunnel.  Moreover, since a resulting
compressing disk for $Q_{x}$ lies in $N$ and so cannot persist into
$Q$, the tunnel attaching $Q_{0}$ must be incident to that same
boundary component of $Q_{1}$.  It is easy to see then that $Q$ is a
tube-spanned recessed product, where the two recessed surfaces are
$Q_{0}$ and the union of $Q_{1}, Q_{2}$ along their parallel boundary
components.

Similar arguments apply if $Q_{-}$ has one or two components.
This completes the proof of the Claim.

\bigskip

Compressing a surface does not affect its boundary, so the theorem
follows immediately from Lemma \ref{lemma:Ddisjoint} and Proposition
\ref{prop:essential} unless the surface $Q_{X}$, obtained by
maximally compressing $Q$ into $X$ has the property that each of its
non-closed components is boundary parallel in $M$.  Of course the
symmetric statement holds also for the surface $Q_{Y}$ obtained by
maximally compressing $Q$ into $Y$; indeed, all the ensuing arguments
would apply symmetrically to $Q_{Y}$ simply by switching labels $X$
and $Y$ throughout.  So henceforth assume that all components of 
$Q_{X}$ are either closed or $\bdd$-parallel.  There are some of the 
latter, since $Q$ has boundary.

Let $Q_{0}$ be an outermost $\bdd$-parallel component of $Q_{X}$ that
is not closed.  That is $Q_{0}$ is a component which
is parallel to a subsurface of $\bdd M$ and no component of $Q_{X}$
lies in the region $R \cong Q_{0} \times I$ of parallelism.  As in the
proof of Lemma \ref{lemma:Ddisjoint}, use the notation $X' \subset X$ and
$Y' \supset Y$ for the two $3$-manifolds into which $Q_{X}$ splits
$M$, noting that, unlike for $Q$, some components of $Q_{X}$ may be
closed.  Note also that $\Ggg \subset Y'$.

{\bf Case 1:} Some such outermost region $R$ lies in $Y'$ 

In this case the other side of $Q_{0}$ lies in $X'$, and so its
interior is disjoint from $\Ggg$.  Since $Q$ is connected, this
implies that all of $Q$ lies in $R$.  In particular, $\Ggg \subset R$,
all compressing disks for $Q$ in $Y$ also lie in $R$, and $Q_{0} =
Q_{X}$.  Let $D \subset M$ be a $\bdd$-reducing disk for $M$ as in
Lemma \ref{lemma:Ddisjoint} so that $\Ggg$ is disjoint from $D$ and $D
\cap Q_{0}$ consists only of arcs that are essential in $Q_{0}$.

 Any outermost such arc in $D$ cuts off a $\bdd$-reducing disk $D_{0}
 \subset D$.  Suppose first that $D_{0}$ lies in $M - R$ and let
 $Q'_{0}$ be the surface created from $Q_{0}$ by $\bdd$-compressing
 along $D_{0}$.  By Lemma \ref{lemma:boundredessential} $Q'_{0}$ is
 incompressible, so all boundary components of $Q'_{0}$ are essential
 in $\bdd M$ unless $Q_{0}$ is an annulus that is parallel to $\bdd M$
 also via $M - R$.  The latter would imply that $Q_{0}$ is a
 longitudinal annulus of a solid torus, $D$ is a meridian of that
 solid torus and we could have taken for $D_{0}$ the half of $D$ that
 does lie in $R$.  In the general case, the union of $D_{0}$ with a
 disk of parallelism in $R$ gives a $\bdd$-reducing disk for $M$ that
 is disjoint from $\bdd Q'_{0}$ so for any boundary component $q'$ of
 $Q'_{0}$, $d(q', \mathcal{V}) \leq 1$.  Then for $q$ any component of
 $\bdd Q = \bdd Q_{X} = \bdd Q_{0}$, $d(q, \mathcal{V} \leq d(q',
 \mathcal{V}) + d(q, q') \leq 2 \leq  1 - \chi(Q)$ and we are done.
 In any case, we may as well then assume that $D_{0}$ lies in $R
 \subset Y'$.
 
Since $\Ggg$ is disjoint from $D_{0}$, $D_{0}$ is a $\bdd$-reducing disk
for $Q$ as well, lying in $Y$.  Then a standard outermost arc
argument in $D_{0}$ shows that a compressing disk for $Q$ in $Y$ can
be disjoint from $D_{0}$.  Then $\bdd$-reducing $Q$ along $D_{0}$
leaves a surface that is still bicompressible (for meridians of $\Ggg$
constitute compressing disks in $X$) but with $1 - \chi(Q)$ reduced.
The proof then follows by induction.  (In fact, this argument can be
enhanced to show directly that Case 1 simply cannot arise.)

It remains to consider the case in which all outermost
components of $Q_{X}$ are $\bdd$-parallel via a region that lies in
$X'$.  We distinguish two further cases:

{\bf Case 2:} There is nesting among the non-closed components of
$Q_{X}$.  We will prove then that $Q$ must be a tube-spanned recessed 
collar.  

In this case, let $Q_{1}$ be a component that is not closed (so it is
$\bdd$-parallel) and is ``second-outermost''.  That is, the region of
parallelism between $Q_{1}$ and $\bdd M$ contains in its interior only
outermost components of $Q_{X}$; denote the union of the latter
components by $Q_{0}$.  Then the region between $Q_{0}$ and $Q_{1}$ is
itself a product $R \cong Q_{1} \times I$ but one end contains $Q_{0}$
as a possibly disconnected subsurface.  Since outermost components cut
off regions lying in $X'$, $R \subset Y'$.  We now argue much as in
Case 1: Since $\Ggg \subset Y'$ and $Q$ is connected, all of $\Ggg$
must lie in $R$, so $Q_{X} = Q_{1} \cup Q_{0}$.  Let $D$ be a
$\bdd$-reducing disk for $M$ that is disjoint from $\Ggg$ and
intersects $Q_{X}$ only in arcs that are essential in $Q_{X}$.  As in
Case 1, each outermost arc of $D \cap Q_{X}$ in $D$ lies in $Q_{0}$.

Choose a complete collection of $\bdd$-compressing disks
$\mathcal{F}$, in the region of parallelism between $Q_{1}$ and $\bdd
   M$, so that the complement $Q_{1}- \mathcal{F}$ is a single disk
$D_{Q}$.  Each disk in $\mathcal{F}$ is incident to $Q_{1}$ in a 
single arc. Now import the argument of Lemma
\ref{lemma:Ddisjoint} into this context: Let $E$ be a compressing disk
for $Y$, here chosen so that $E \cap \mathcal{F}$ is minimized.  This means
first of all that $E \cap \mathcal{F}$ is a collection of arcs. 
As in the proof of Lemma \ref{lemma:Ddisjoint}, $\Ggg$ may be slid and
isotoped so it is disjoint from $\mathcal{F}$.  $\Ggg$ is incident to
$Q_{1}$ since $Q$ is connected.  Since $D_{Q}$ is connected, the ends
of $\Ggg$ on $D_{Q}$ may be slid within $D_{Q}$ so that ultimately
$\Ggg$ is incident to $D_{Q}$ in a single point.  $\bdd E$ is
necessarily incident to that end, since $Q$ is weakly incompressible.
It follows that $\bdd E$ cannot be incident to $Q$ only in $D_{Q}$
(else $\bdd E$ could be pushed off the end of $\Ggg$ in $D_{Q}$) so
$\bdd E$ must intersect the arcs $\bdd \mathcal{F} \cap Q_{1}$.  Let
$\bbb \subset (\mathcal{F} \cap E)$ be outermost in $E$ among all arcs
incident to components of $\bdd \mathcal{F} \cap Q_{1}$.  Let $E_{0}$
be the disk that $\bbb$ cuts off from $E$.

If both ends of $\bbb$ were in $\mathcal{F} \cap Q_{1}$ then, since
each disk of $\mathcal{F}$ is incident to $Q_{1}$ in a single arc,
$\bbb$ would cut off a subdisk of $\mathcal{F}$ that could be used to
alter $E$, creating a compressing disk for $Y$ that intersects
$\mathcal{F}$ in fewer points.  We conclude that the other end of
$\bbb$ is on $Q_{0}$.  Since $\bbb$ is outermost among those arcs of
$E \cap \mathcal{F}$ incident to $D_{Q}$, $\bdd E_{0}$ traverses the
end of $\Ggg$ on $D_{Q}$ exactly once.  So, as in the proof of Lemma
\ref{lemma:Ddisjoint}, it can be used to slide and isotope an edge $\rho$ of
$\Ggg$ until it coincides with $\bbb$.  Hence the edge $\rho \subset
\Ggg$ can be made into a vertical arc (i.  e. an $I$-fiber) in the
product structure $R = Q_{1} \times I$.

Using that product structure and an essential circle in the component
of $Q_{0}$ that is incident to $\rho$, $\rho$ can be viewed as part of
a vertical incompressible annulus $A$ with ends on $Q_{1}$ and
$Q_{0}$.  Now apply the argument of Lemma \ref{lemma:Ddisjoint} again:
$A - \rho$ is a disk $E'$.  Since $E'$ is a disk, use the argument of Lemma
\ref{lemma:Ddisjoint} to slide and isotope the edges of $\Ggg - \rho$
until they are disjoint from $E'$.  After these slides, $E'$ is revealed
as a compressing disk for $Q$ in $Y$.  On the other hand, if there is
in fact any edge $\ggg$ in $\Ggg - \rho$, the compressing disk for $Q$
in $X$ given by the meridian of $\eta(\ggg)$ would be disjoint from
$E$, contradicting weak incompressiblilty of $Q$.  So we conclude that
in fact $\Ggg = \rho$ and so, other than the components of $Q_{X}$
incident to the ends of $\rho$, each component of $Q_{X}$ is a
component of $Q$; since $Q$ is connected, there are no such other 
components.  That is, $Q$ is obtained by tubing $Q_{1}$ to
the connected $Q_{0}$ along $\rrr$ and so is a tube-spanned recessed
collar.  This completes the argument
in this case.

{\bf Case 3:} All non-closed components of $Q_{X}$ are outermost among
the components of $Q_{X}$.  We will show that in this case $Q$ is
strongly $\bdd$-compressible; the proof then follows from the Claim
above.

We have already seen that all non-closed components of $Q_{X}$ are
$\bdd$-parallel through $X'$.  Choose a $\bdd$-reducing disk $D
\subset M$ as in Lemma \ref{lemma:Ddisjoint} so that $D$ is disjoint
from the graph $\Ggg$, intersects $Q_{X}$ mimimally and intersects $Q$
only in arcs that are essential in $Q_{X}$.  Although there is no
nesting among the components of $Q_{X}$, it is not immediately clear
that the arcs $D \cap Q_{X}$ are not nested in $D$.  However, it is
true that each outermost arc cuts off a subdisk of $D$ that lies in
$X'$, as shown in the proof of Case 1 above.  In what follows, $D'$
will represent either $D$, if no arcs of $D \cap Q_{X}$ are nested in
$D$, or a disk cut off by a ``second-outermost'' arc of intersection
$\lll_{0}$ if there is nesting.  Let $\Lll \subset D'$ denote the
collection of arcs $D' \cap Q$; one of these arcs (namely $\lll_{0}$)
may be on $\bdd D'$.

\begin{figure}[tbh]
\centering
\includegraphics[width=0.6\textwidth]{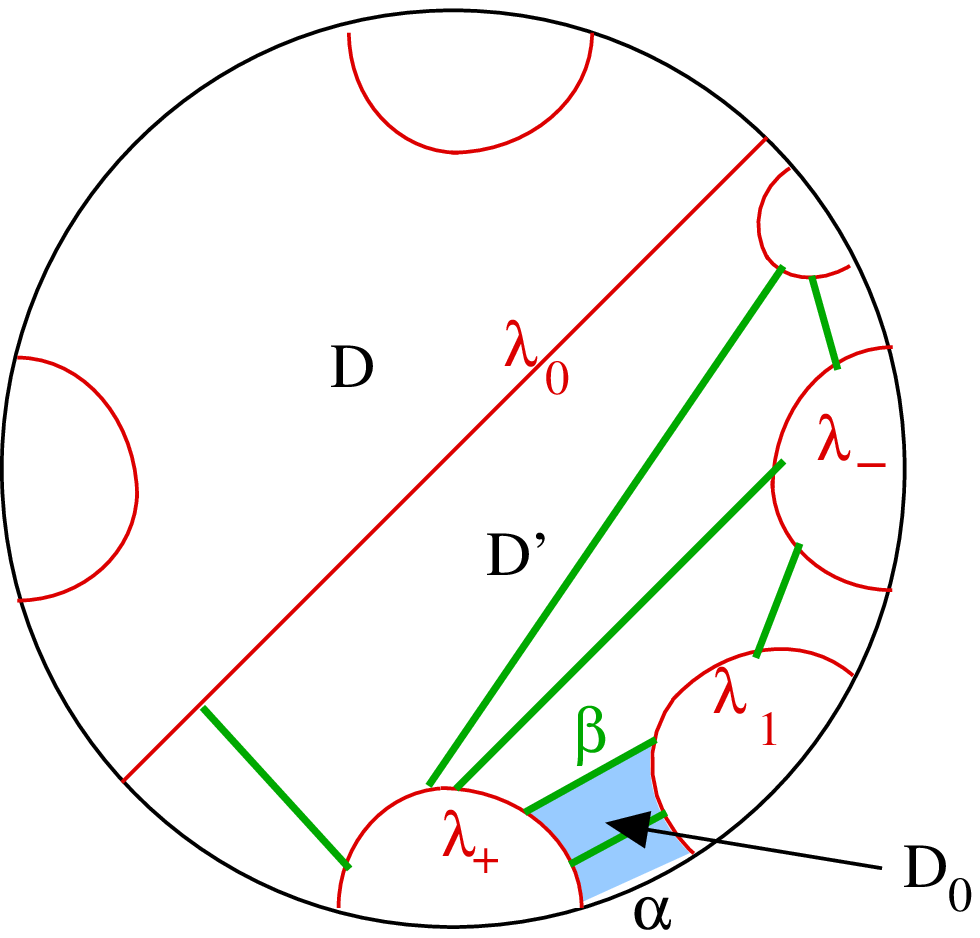}
\caption{} \label{fig:rectangle}
\end{figure}

Consider how a compressing disk $E$ for $Q$ in $Y$ intersects $D'$.
All closed curves in $D' \cap E$ can be removed by a standard
innermost disk argument redefining $E$.  Any arc in $D' \cap E$ must
have its ends on $\Lll$; a standard outermost arc argument can be used
to remove any that have both ends on the same component of $\Lll$.  If
any component of $\Lll - \lll_{0}$ is disjoint from all the arcs $D'
\cap E$, then $Q$ could be $\bdd$-compressed without affecting $E$.
This reduces $1 - \chi(Q)$ without affecting bicompressibility, so we
would be done by induction.  Hence we restrict to the case in which
each arc component of $\Lll - \lll_{0}$ is incident to some arc
components of $D' \cap E$.  See Figure \ref{fig:rectangle}.

It follows that there is at least one component $\lll_{1} \neq \lll_{0}$ of
$\Lll$ with this property: any arc of $D' \cap E$ that has one end
incident to $\lll_{1}$ has its other end incident to one of the (at most
two) neighboring components $\lll_{\pm}$ of $\Lll$ along $\bdd D'$.
(Possibly one or both of $\lll_{\pm}$ are $\lll_{0}$.)  Let $\bbb$ be the
outermost arc in $E$ among all arcs of $D' \cap E$ that are incident to
the special arc $\lll_{1}$.  We then know that the other end of $\bbb$ is
incident to (say) $\lll_{+}$ and that the disk $E_{0} \subset E$ cut
off by $\bbb$ from $E$, although it may be incident to $D$ in its
interior, at least no arc of intersection $D \cap interior(E_{0})$ is
incident to $\lll_{1}$.

Let $D_{0}$ be the rectangle in $D$ whose sides consist of subarcs of
$\lll_{1}$, $\lll_{+}$, $\bdd D$ and all of $\bbb$.  Although $E$ may
intersect this rectangle, our choice of $\bbb$ as outermost among arcs
of $D \cap E$ incident to $\lll_{1}$ guarantees that $E_{0}$ is
disjoint from the interior of $D_{0}$ and so is incident to it only in
the arc $\bbb$.  The union of $E_{0}, D_{0}$ along $\bbb$ is a disk
$D_{1} \subset Y$ whose boundary consists of the arc $\aaa = \bdd M
\cap \bdd D_{0}$ and an arc $\bbb' \subset Q$.  The latter arc is the
union of the two arcs $D_{0} \cap Q$ and the arc $E_{0} \cap Q$.  If
$\bbb'$ is essential in $Q$, then $D_{1}$ is a $\bdd$-compressing disk
for $Q$ in $Y$ that is disjoint from the boundary compressing disk in
$X$ cut off by $\lll_{1}$.  So if $\bbb'$ is essential then $Q$ is
strongly $\bdd$-compressible and we are done by the Claim.

Suppose finally that $\bbb'$ is inessential in $Q$.  Then $\bbb'$ is
parallel to an arc on $\bdd Q$ and so, via this parallelism, the disk
$D_{1}$ is itself parallel to a disk $D'$ that is disjoint from $Q$
and either is $\bdd$-parallel in $M$ or is itself a $\bdd$-reducing
disk for $M$.  If $D'$ is a $\bdd$-reducing disk for $M$, then $\bdd 
D' \in \mathcal{V}$,  $d(\mathcal{Q}, \mathcal{V})
\leq 1$ and we are done.  On the other hand, if $D'$ is parallel to
a subdisk of $\bdd M$, then an outermost arc of $\bdd D$ in that disk
(possibly the arc $\aaa$ itself) can be removed by an isotopy of $\bdd
D$, lowering $|D \cap Q| = |D \cap Q_{X}|$.  This contradiction to our
original choice of $D$ completes the proof.  \qed

\end{document}